\documentclass[12pt]{article}
\usepackage{amsmath}
\usepackage{amssymb}

\usepackage[cp1251]{inputenc}
\usepackage[russian]{babel}
\begin{document}
\vskip 12pt\noindent
MSC 34C10, 34C11, 34C99
\vskip 20pt
{\bf\centerline {Stability criteria for   second order linear}
\centerline {ordinary differential equations}
\vskip 20pt
\centerline { G. A. Grigorian}}
\vskip 20pt
Abstract. We use some properties of solutions of Riccati equation for establishing   boundedness and stability criteria for  solutions of   second order linear ordinary differential equations. We show that the conditions on coefficients of the equations, appearing in the proven criteria, do not follow from the conditions, which  ensure the application of the WKB approximation to the second order linear equations. On these examples we  compare the obtained results wit the  results obtained by the Liapunov and Bogdanov methods, by a method involving estimates of solutions in the Lozinski's logarithmic norms, and by the freezing method. We compare these results with the Wazevski's theorem as well.

\vskip 20pt

\vskip 20pt
{\bf \centerline { \S 1. Introduction}}

\vskip 20pt

Let $p(t)$  and  $q(t)$  be complex valued continuous functions on $[t_0;+\infty)$. Consider the equation
$$
\phi''(t) + p(t)\phi'(t) + q(t)\phi(t) = 0,\phantom{aaa} t\ge t_0. \eqno (1.1)
$$
Study  of the boundedness and stability behavior of  solutions of Eq.  (1.1)  is an important problem of the qualitative theory of differential equations and many works are devoted to it (see e.g., the book [1]  and cited works therein, [2  - 12]).

Let  $p(t)$ be continuously differentiable. In Eq. (1.1)  make the substitution
$$
\phi(t) = E(t)\psi(t), \phantom{aaa} t\ge t_0, \eqno (1.2)
$$
where  $E(t)\equiv \exp\biggl\{-\frac{1}{2}\int\limits_{t_0}^t p(\tau)d\tau\biggr\}$.
We get
$$
\psi''(t) - \frac{D(t)}{4}\psi(t)= 0, \phantom{aaa}t\ge t_0, \eqno  (1.3)
$$
where $D(t)\equiv 2p'(t) + p^2(t) - 4 q(t),\phantom{a} t\ge t_0$.  One of important methods of studying the  boundedness and stability problems of the solutions of Eq. (1.1) is the application  of the Liouville's  transformation   (see [2], pp.  131, 132, 152, 153). In the book  [5] on the basis of the Liouville's transformation a  substantiation  of asymptotic representation of the solutions  of Eq. (1.3) and their derivatives  is given (see. [5], pp. 54 - 61, WKB approximation   [Wentzel–Kramers–Brillouin]). It is assumed therein, that $D(t)$  is twice continuously differentiable, $D(t)\neq 0,\phantom{a} t\ge t_0,  Re \sqrt{D(t)} \ge 0$  for  $t>>1$   and
$$
\int\limits_{t_0}^{+\infty}\left|36 \frac{D''(\tau)}{D(\tau)^{3\slash 2}} - 5 \frac{D'(\tau)^2}{D(\tau)^{5\slash 2}}\right|d\tau < +\infty. \eqno (1.4)
$$
By virtue of (1.2) the WKB  approximation gives possibility to describe wide classes of equations (1.1) with bounded and (or) unbounded solutions, classes of stable and (or) unstable equations (1.1) in terms of their coefficients.

Assume $x(t)$ is a nonnegative continuous function on the half line $[t_0;+\infty)$. Consider the Riccati equation
$$
y'(t) + y^2(t) = x(t), \phantom{aaa} t\ge t_0. \eqno (1.5)
$$

{\bf Definition 1.1}. {\it The solution $y(t)$ of Eq. (1.5) satisfying the initial condition\linebreak $y(t_0) = \sqrt{x(t_0)}$ is called differential root.}

For the study of the boundedness and stability problem of solutions of Eq. (1.1) in this work the Riccati equations method is applied, which (in this work) basically  is an application of properties of the differential root of $\frac{D(t)}{4}$, corresponding to the solutions of Eq. (1.3). Unlike conditions on $D(t)$, providing of use WKB approximation, here other restrictions are imposed on $D(t)$ assuming the condition

\noindent
A) $D(t) > 0, \phantom{a} t\ge t_0,$ and $p(t), \phantom{a} D(t)$ are  continuously differentiable functions;

\noindent
and other conditions, different from (1.4), be satisfied.
Note that the case  $D(t)<0, \phantom{a} t\ge t_0$, is studied in [10].
Boundedness and stability tests for the solutions of Eq. (1.1) in terms of their coefficients are proved. Examples, to which the mentioned tests are applicable and which do not satisfy the condition (1.4), are represented.

\vskip 20pt
{\bf \centerline { \S 2. Main results}}

\vskip 20pt

For any positive and continuously differentiable on $[t_0;+\infty)$ function $x(t)$ denote

$$
R_x(t_1;t)\equiv \frac{1 + \sqrt{x(t_0)}(t_1 - t_0)}{1 + \sqrt{x(t_0)}(t - t_0)}\exp\biggl\{-\int\limits_{t_1}^t\sqrt{x(s)} ds\biggr\}\sup\limits_{\xi \in [t_0;t_1]}\frac{|(\sqrt{x(\xi)})'|}{\sqrt{x(\xi)}}
+\sup\limits_{\xi \in [t_1;t]}\frac{|(\sqrt{x(\xi)})'|}{\sqrt{x(\xi)}},
$$
$$
\rho_x(t) \equiv \inf\limits_{t_1\in [t_0;t]}R_x(t_1;t),\phantom{aaa}  t_0 \le t_1 \le t,
$$

In our main results  the functions

$$
\rho_{_{D\slash4}}(t), \phantom{aaa}    r_1(t)\equiv \int\limits_{t_0}^t [\sqrt{D(\tau)} - Re\hskip 2pt p(\tau)] d\tau - \frac{1}{2} \ln D(t),
$$
$$
r_2(t)\equiv \int\limits_{t_0}^t [\sqrt{D(\tau)} - Re\hskip 2pt p(\tau)] d\tau - \frac{1}{2} \ln D(t) + 2 \ln [1 + |p(t) - \sqrt{D(t)}|],\phantom{aaa} t\ge t_0;
$$
play a crucial role

{\bf Theorem 2.1}. {\it Let  the conditions

\noindent
A) $D(t) > 0, \phantom{a} t\ge t_0,$ and $p(t), \phantom{a} D(t)$ are  continuously differentiable functions;

\noindent
and one of the following groups of conditions

\noindent
B) $D(t)$  is a nondecreasing function;  for some $\varepsilon > 0$
the function $\frac{D'(t)}{D(t)^{3\slash 2 - \varepsilon}}$  is bounded;

\noindent
C) $D(t)\ge \varepsilon > 0,\phantom{a} t\ge t_0$,   the function $\frac{D'(t)}{D(t)}$  is bounded  and $\int\limits_{t_0}^{+\infty}\rho_{_{D\slash4}}(s) \frac{|D'(s)|}{D(s)^{3\slash 2}} d s < +\infty$

\noindent
be satisfied.
Then  all solutions of Eq. (1.1) are bounded (vanish on $+\infty$) if and only if
the function $r_1(t)$ is bounded from above ($\lim\limits_{t\to +\infty} r_1(t) = - \infty$).}

In many cases in applications of Eq. (1.1)  its stability property plays an important role, and the property of boundedness of its solutions is a necessary condition for stability of Eq. (1.1). However this property (even the property of vanishing of all solutions to Eq. (1.1) in $+\infty$) still does  not  guarantee the stability of Eq. (1.1). The next theorem indicates some conditions on the coefficients of Eq. (1.1) which guarantee  Liapunov stability (asymptotically stability) of Eq. (1.1).

{\bf Theorem 2.2.} {\it Let the condition A) and  the group of conditions C) of  Theorem~2.1 or the  group of conditions

\noindent
D) $D(t)$ is a nondecreasing function;
$\frac{D'(t)}{D(t)}$   is bounded,

\noindent
be satisfied. Then  Eq. (1.1) is  Liapunov  stable (asymptotically) if and only if  the function $r_2(t)$  is bounded from above ($\lim\limits_{t\to+\infty} r_2(t)= - \infty$).}


{\bf Corollary 2.1.} {\it Let  $D(t)\ge \varepsilon > 0,\phantom{a} t\ge t_0; \frac{|D'(t)|}{D(t)} \le \frac{c}{(1 + t - t_0)^\alpha}, \phantom{a}t\ge t_0,\phantom{a} c>0,\linebreak  \alpha>0$;
  \phantom{a} $\int\limits_{t_0}^{+\infty}\frac{d\tau}{\sqrt{D(\tau)}(1 + \tau - t_0)^{2\alpha}} < +\infty$ and let the condition A) be satisfied.  Then the following  assertions are valid:

\noindent
А$_1$) All solutions of Eq. (1.1) are bounded
(vanish on $+\infty$ ) if and only if
the function  $r_1(t)$  is bounded  above ($\lim\limits_{t\to+\infty} r_1(t)= - \infty$).

\noindent
B$_1$)  Eq. (1.1) is Liapunov stable (asymptotically) if and only if  the function $r_2(t)$ is bounded  above ($\lim\limits_{t\to+\infty} r_2(t)= - \infty$).}

Example 2.1. Consider the equation
$$
\phi''(t) + p_1(t) \phi'(t) + q_1(t) \phi(t) = 0, \phantom{aaa} t\ge 1, \eqno (2.1)
$$
where $p_1(t) \equiv \lambda t, \phantom{a} q_1(t) \equiv \frac{\lambda}{2} + \frac{\lambda^2 t^2}{4} - \frac{t}{4} - \frac{1}{4} \int\limits_1^t\sin^2 e^\tau d\tau, \phantom{a} t\ge 1, \phantom{a} \lambda = const \in \mathrm{\bf C}$. For this equation we have $ D(t) = D_1(t)\equiv t+ \int\limits_1^t\sin^2 e^\tau d \tau,\phantom{a} t\ge 1$, and $D_1(t)$ is an increasing function on $[1;+\infty)$. Therefore for Eq. (2.1) the conditions A) and B) hold. For Eq. (2.1) we have
$$
r_1(t) = \int\limits_1^t\biggl[\sqrt{\tau + \int _1^\tau \sin^2 e^s d s} - Re \lambda \tau\biggr] d \tau - \frac{1}{2} \ln \biggl[ t + \int_1^t \sin^2 e^\tau d\tau\biggr], \phantom{aaa} t\ge 1.
$$
$$
\mbox{Hence}\phantom{a} \lim\limits_{t \to +\infty} r_1(t) =\left\{
\begin{array}{l}
- \infty,\phantom{a} if \phantom{a} Re \lambda > 0;\\
+\infty,\phantom{a} if\phantom{a} Re \lambda \le 0.
\end{array}\right.\phantom{a} \mbox{Therefore by Theorem 2.1 if}\phantom{a} Re \lambda > 0\phantom{a}\mbox{then}
$$
all solutions of Eq. (2.1) vanish on $+\infty$ and if $Re \lambda \le 0$ then Eq. (2.1) has unbounded solution, i.e., Eq. (2.1) is unstable. It is not difficult to show, that for  $D(t)=D_1(t)$ the condition (1.4) does not hold. Therefore the WKB approximation is not applicable to Eq. (2.1). The substitution $\phi'(t) = \psi(t), \phantom{a} t\ge 1$, in Eq. (2.1) reduces it to the system
$$
\left\{
\begin{array}{l}
\phi'(t) = \phantom{aaaaaaaaaaaaaaa} \psi(t);\\
\psi'(t) = - q_1(t) \phi(t) - p_1(t) \psi(t), \phantom{a} t\ge 1.
\end{array}
\right.
$$
It is not difficult to verify that the application of estimates of Liapunov ([4], p. 132) and Bogdanov ([4], p. 133), the estimate by Lozinski's logarithmic norms ([4], p. 137), as well as the estimation by freezing method ([4], p. 139) to the last system give no result. The application of Wazevski's theorem to the last system also gives no result. Hence these estimates and the Wazvski's theorem give no rezult for Eq. (2.1).

Example 2.2. Consider the equation
$$
\phi''(t) + p_2(t) \phi'(t) + q_2(t) \phi(t) = 0, \phantom{aaa} t\ge 1, \eqno (2.2)
$$
where $p_2(t) \equiv \lambda t^2, \phantom{a} q_2(t) \equiv \lambda t + \frac{\lambda^2 t^4}{4} - \frac{t^2}{4} - \frac{1}{4} \biggl(\int\limits_1^t \sin e^\tau d\tau\biggr)^2, \phantom{a} t\ge1, \phantom{a} \lambda = const \in \mathrm{\bf C}$.
For this equation we have $D(t)=D_2(t)=t^2 +  \left(\int\limits_1^t\sin e^\tau d \tau\right)^2, \phantom{a} r_2(t) = \int\limits_1^t\Bigl[\sqrt{\tau^2 + \bigl(\int_1^\tau \sin e^s d s\bigr)^2} - -Re \lambda \tau^2\Bigr] d \tau - \frac{1}{2} \ln D_2(t) + 2\ln[1+ |p_2(t) - \sqrt{D_2(t)}|], \phantom{a}t\ge 1$.
It is not difficult to check that the conditions A) and C) for Eq. (2.2) hold and for $D(t) =D_2(t)$ the condition (1.4) does not fulfill. Therefore Theorem 2.2 is applicable to Eq. (2.2) and  the WKB approximation is not applicable to Eq. (2.2). We have
$$
\lim\limits_{t \to +\infty} r_2(t) =\left\{
\begin{array}{l}
- \infty,\phantom{a} if \phantom{a} Re \lambda > 0;\\
+\infty,\phantom{a} if\phantom{a} Re \lambda \le 0.
\end{array}\right.
$$
By Theorem 2.2 from here it follows that for  $Re \lambda > 0$ Eq. (2.2) is asymptotically stable and for
$Re \lambda \le 0$ Eq. (2.2) is unstable. Moreover we also can use Theorem 2.1 to Eq. (2.2) and show that for  $Re \lambda \le 0$ it has an unbounded solution. It is not difficult to verify that the application of the mentioned above estimates and the Wazevski's theorem to  Eq. (2.2) gives no result.

Example 2.3. Consider the equation
$$
\phi''(t) + p_3(t) \phi'(t) + q_3(t) \phi (t) = 0, \phantom{aaa} t\ge 1, \eqno (2.3)
$$
where $p_3(t) \equiv \lambda + \mu \sin t, \phantom{a} q_3(r) \equiv \frac{\mu \cos t}{2} + \frac{(\lambda + \mu \sin t)^2}{4} - \frac{1}{4}\biggl(\alpha + \beta \cos\ln t + \gamma \int\limits_1^t \frac{\sin^2\tau}{\tau} d\tau\biggr), \linebreak \lambda = const \in \mathrm{\bf C},\phantom{a} \mu = const \in \mathrm{\bf C},\phantom{a} \alpha = const \ge \beta = const > 0,\phantom{a} \gamma = const >~ 0.$ For this equation we have $D(t) = D_3(t) = \alpha + \beta \cos\ln t + \gamma \int_1^t \frac{\sin^2 \tau}{\tau} d \tau, \linebreak r_1(t) = \int\limits_1^t \biggl[\sqrt{\alpha + \beta \cos\ln \tau + \gamma \int_1^\tau \frac{\sin^2 s}{s} d s} -  Re \lambda - Re \mu \sin\tau \biggr]d\tau - \frac{1}{2}\ln \biggl(\alpha + \beta \cos\ln t + \gamma \int\limits_1^t \frac{\sin^2\tau}{\tau} d\tau \biggr), \phantom{a} r_2(t) = r_1(t) + 2 \ln\biggl[1 + \biggl|\lambda + \mu \sin t - \alpha - \beta \cos\ln t - \gamma \int\limits_1^t \frac{\sin^2\tau}{\tau} d\tau\biggr|\biggr], \phantom{a} t\ge 1.$ Hence
$$
\lim\limits_{t\to +\infty} r_1(t) = \lim\limits_{t\to +\infty} r_2(t) = \left\{
\begin{array}{l}
- \infty,\phantom{a} if \phantom{a} Re \lambda > \sqrt{\alpha} ;\\
+\infty,\phantom{a} if\phantom{a} Re \lambda \le \sqrt{\alpha}.
\end{array}\right. \eqno (2.4)
$$
We can easily check that for $D(t) = D_3(t)$ condition (1.4) does not hold. Therefore the WKB approximation is not applicable to Eq. (2.3). It is not difficult to verify that for Eq. (2.3) all conditions of Corollary 2,1 are fulfilled. Therefore taking into account (2.4) we  get:

\noindent
for $Re \lambda > \sqrt{\alpha}$ Eq. (2.3) is asymptotically stable;

\noindent
for $Re \lambda \le \sqrt{\alpha}$ Eq. (2.3) has unbounded solution.

It is not difficult to verify that the application of the mentioned above estimates and the Wazevski's theorem to  Eq. (2.3) gives no result. Note that the results of work [11] concern to the case  $D(t) < 0, \phantom{a} t\ge t_0,$ and the results of work [12] concern to the case of periodic functions $p(t)$ and $q(t)$. Therefore the results of these works cannot be applicable to the equations (2.1) - (2.3).

\vskip 20pt

{\bf \centerline { \S 3. Proof of the main results}}

\vskip 20pt

To prove the main results at fist we shall formulate and prove some preliminary propositions.
Let  $x_1(t)$  be a real valued continuous  function on $[t_0;+\infty)$ . Along with Eq. (1.5)  consider the Riccaty equation
$$
y'(t) + y^2(t) = x_1(t),\phantom{aaa} t\ge t_0.  \eqno (3.1)
$$
The following assertion is valid (see [13]).

{\bf Theorem 3.1}. {\it Let Eq. (1.5) has a real valued solution  $y_0(t)$    on  $[t_0;+\infty)$, and let  $x_1(t) \ge x(t),\phantom{a} t\ge~t_0$.  Then for each  $y_{(0)} \ge y_0(t_0)$
Eq. (3.1) has a solution  $y_1(t)$  on $[t_0;+\infty)$,
satisfying the initial  condition $y_1(t_0)=y_{(0)}$, moreover  $y_1(t) \ge y_0(t),\phantom{a} t\ge t_0$.}

The proof of a more general theorem is presented in [14].

Since  $y_0(t)\equiv 0$   is a solution of the equation
$$
y'(t) + y^2(t) = 0,\phantom{aaa} t\ge t_0,
$$
from  Theorem 3.1 we immediately get:

{\bf Corollary 3.1}.  {\it Let  $ x(t) \ge 0, \phantom{a} t\ge t_0$.
Then for any $y_{(0)} \ge 0$  Eq.  (1.5)
has a \linebreak solution $y_1(t)$ on  $[t_0;+\infty)$,
satisfying the initial  condition  $y_1(t_0) = y_{(0)}$, moreover \linebreak  $y_1(t) \ge 0,\phantom{a} t\ge  t_0$.}

From Corollary 3.1 it follows, that the differential root is defined on  $[t_0;+\infty)$  and is nonnegative.

{\bf Remark 3.1} {\it A more detailed study of the properties of the differential root is presented in [13].}

In the sequel the differential root of  $x(t)$ we shall denote by  $y_x(t)$.

Let  $x(t)$ be continuously differentiable and  $x(t) > 0, \phantom{a}t\ge t_0 $. Then
$$
[y_x(t) - \sqrt{x(t)}]' + (y_x(t) + \sqrt{x(t)})[y_x(t) - \sqrt{x(t)}] = - (\sqrt{x(t)})',\phantom{aaa} t\ge t_0.
$$
It follows from here, that  $u_0(t) \equiv y_x(t) - \sqrt{x(t)}\phantom{a} (t\ge t_0)$
is a solution of the first order linear equation:
$$
u'(t) + F(t) u(t) = - (\sqrt{x(t)})', \phantom{aaa} t\ge t_0,
$$
where  $F(t) \equiv y_x(t) + \sqrt{x(t)},\phantom{a} t\ge t_0$. Therefore by Cauchy formula
$$
y_x(t) - \sqrt{x(t)} =exp\biggl\{- \int\limits_{t_1}^tF(\tau) d \tau\biggr\}\times \phantom{aaaaaaaaaaaaaaaaaaaaaaaaaaaaaaaaaaaaaaaaaaaaaa}
$$
$$
 \times\left[y_x(t_1) - \sqrt{x(t_1)} - \int\limits_{t_1}^t\exp\biggl\{\int\limits_{t_1}^\tau F(s) d s\biggr\}\bigl(\sqrt{x(\tau)}\bigr)'d \tau\right],\phantom{aaa} t,\phantom{a} t_1 \ge t_0, \eqno (3.2)
$$
in particular,
$$
y_x(t) - \sqrt{x(t)} = - \int\limits_{t_0}^t\exp\biggl\{-\int\limits_\tau^t F(s) d s\biggr\}\bigl(\sqrt{x(\tau)}\bigr)'d \tau,\phantom{aaa} t\ge t_0.
$$
Hence 
$$
|y_x(t_1) - \sqrt{x(t_1)}|  = \int\limits_{t_0}^{t_1} F(\tau)\exp\biggl\{- \int\limits_\tau ^{t_1} F(s) d s\biggr\} \frac{(\sqrt{x(\tau)})'}{F(\tau)} d \tau \le  \sup\limits_{\xi\in [t_0;t_1]}\frac{|(\sqrt{x(\xi)})'|}{F(\xi)} \times \phantom{aaaaa}
$$
$$
\times \int\limits_{t_0}^{t_1} d\biggl[ \exp\biggl\{- \int\limits_\tau ^{t_1} F(s) d s\biggr\}\biggr]  = \sup\limits_{\xi\in [t_0;t_1]}\frac{|(\sqrt{x(\xi)})'|}{\sqrt{x(\xi)}} \biggl[1 - \exp\biggl\{- \int\limits_{t_0} ^{t_1} F(s) d s\biggr\}\le
 $$
 $$
\phantom{aaaaaaaaaaaaaaaaaaaaaaaaaaaaaaaaaaaaaaaaaaaaaaaaaa}\le \sup\limits_{\xi\in [t_0;t_1]}\frac{|(\sqrt{x(\xi)})'|}{\sqrt{x(\xi)}}.
$$
Then from (3.2) we get:
$$
|y_x(t) - \sqrt{x(t)}| \le \phantom{aaaaaaaaaaaaaaaaaaaaaaaaaaaaaaaaaaaaaaaaaaaaaaaaaaaaaaa}
$$
$$
\phantom{aaaaaa}\le \exp\biggl\{- \int\limits_{t_1}^t F(\tau) d\tau\biggr\}\sup\limits_{\xi \in [t_0;t_1]}\frac{|(\sqrt{x(\xi)})'|}{\sqrt{x(\xi)}} +
\int\limits_{t_1}^t\exp\biggl\{- \int\limits_\tau ^t  F(s)d s\biggr\}|(\sqrt{x(\tau)})'| d\tau\le
$$
$$
\le \exp\biggl\{- \int\limits_{t_1}^t F(\tau) d\tau\biggr\}\sup\limits_{\xi \in [t_0;t_1]}\frac{|(\sqrt{x(\xi)})'|}{\sqrt{x(\xi)}} +
\sup\limits_{\xi \in [t_1;t]}\frac{|(\sqrt{x(\xi)})'|}{\sqrt{x(\xi)}},    \eqno (3.3)
$$
$t_0\le t_1 \le t,$ as far as
$$
\int\limits_{t_1}^t\exp\biggl\{- \int\limits_\tau ^t  F(s)d s\biggr\}|(\sqrt{x(\tau)})'| d\tau = \int\limits_{t_1}^t F(\tau)\exp\biggl\{- \int\limits_\tau ^t  F(s)d s\biggr\}\frac{|(\sqrt{x(\tau)})'|}{F(\tau)} d\tau \le
$$
$$
\le \sup\limits_{\xi \in [t_0;t_1]}\frac{|(\sqrt{x(\xi)})'|}{F(\xi)}\int\limits_{t_1}^t F(\tau)\exp\biggl\{- \int\limits_\tau ^t  F(s)d s\biggr\}d\tau =
$$
$$
= \sup\limits_{\xi \in [t_0;t_1]}\frac{|(\sqrt{x(\xi)})'|}{F(\xi)}\int\limits_{t_1}^t d\biggl[\exp\biggl\{- \int\limits_\tau ^t  F(s)d s\biggr\}\biggr] = \sup\limits_{\xi \in [t_0;t_1]}\frac{|(\sqrt{x(\xi)})'|}{F(\xi)}\biggl[1 -
$$
$$
\phantom{aaaaaaaaaaaaaaaaaaaaaaaaaaaaaaaaaaa} - \exp\biggl\{- \int\limits_{t_1} ^t  F(s)d s\biggr\}\biggr]
\le \sup\limits_{\xi \in [t_1;t]}\frac{|(\sqrt{x(\xi)})'|}{\sqrt{x(\xi)}}
$$

 {\bf Lemma  3.1}. {\it For every   $s \ge t_0$   the  inequality
$$
y_x(s) \ge \frac{y_x(t_0)}{1 + y_x(t_0)(s - t_0)}.
$$
is valid}

See the proof in  [13].

By virtue of this lemma we have:
$$
\int\limits_{t_1}^t y_x(s) ds \ge \int\limits_{t_1}^t \frac{\sqrt{x(t_0)} d s}{1 + \sqrt{x(t_0)}(s - t_0)}= \ln \frac{1 +\sqrt{x(t_0)} (t - t_0)}{1 + \sqrt{x(t_0)}(t_1 - t_0)},\phantom{aaa} t, t_1 \ge t_0,\phantom{a} t_1 \le t.
$$
From here and from (3.3) it follows:
$$
|y_x(t) - \sqrt{x(t)}| \le \left[\frac{1 + \sqrt{x(t_0)}(t_1 - t_0)}{1 + \sqrt{x(t_0)}(t - t_0)}\exp\biggl\{-\int\limits_{t_1}^t\sqrt{x(s)} ds\biggr\}\sup\limits_{\xi \in [t_0;t_1]}\frac{|(\sqrt{x(\xi)})'|}{\sqrt{x(\xi)}}  +
\right.
 $$
 $$
 \left.+\sup\limits_{\xi \in [t_1;t]}\frac{|(\sqrt{x(\xi)})'|}{\sqrt{x(\xi)}}\right]=R_x(t_1;t),\phantom{aaa} t_0 \le t_1 \le t.
$$
It means that,
$$
|y_x(t) - \sqrt{x(t)}| \le \inf\limits_{t_1\in [t_0;t]}R_x(t_1;t)= \rho_x(t), \phantom{aaa}t\ge t_0.   \eqno (3.4)
$$
If $\frac{1}{2}\frac{|x'(t)|}{x(t)} \le c,\phantom{a} t\ge t_0$,  then it is evident, that
$$
\rho_x(t) \le R_x(t_1;t) \le c,\phantom{aaa} t\ge t_0.      \eqno (3.5)
$$
Let
$$
x(t) \ge \varepsilon > 0,  \frac{|x'(t)|}{x(t)} \le \frac{c}{[1 + \sqrt{x(t_0)}(t - t_0)]^\alpha},\phantom{aaa} t\ge t_0,\phantom{a} c > 0, \phantom{a}\alpha > 0,  \eqno (3.6)
$$
Let us define  $t_1=t_1(t)$  by relation
$$
t - t_1 = \frac{\alpha}{\sqrt{\varepsilon}} \ln[1 + \sqrt{x(t_0)}(t - t_0)],\phantom{aaa} t > \overline{t},
$$
where  $ \overline{t}\phantom{a} (< +\infty)$  satisfies the condition:

\noindent
$ \ln[1 + \sqrt{x(t_0)}(t - t_0)] < \frac{1}{2}(t - t_0), \phantom{a}t > \overline{t}$
(since  $\frac{\ln[1 + \sqrt{x(t_0)}(t - t_0)]}{t - t_0} \to 0$  for  $t \to +\infty$,
the number  $\overline{t}$  always exists). Since  $x(t) \ge \varepsilon > 0,\phantom{a} t\ge t_0$, we have
$$
\int\limits_{t_1}^t\sqrt{x(s)}d s \ge\sqrt{\varepsilon}(t - t_1) = \ln[1 + \sqrt{x(t_0)(t - t_0)}]^\alpha,\phantom{aaa} t_0 \le t_1 \le t.
$$
Therefore, taking into account (3.6) we  get:
$$
R_x(t_1(t);t)\le \frac{1}{2}\frac{1 + \sqrt{x(t_0)}(t_1 - t_0)}{[1 + \sqrt{x(t_0)}(t - t_0)]^{1+\alpha}}\sup\limits_{\xi \in [t_0;t_1]}\frac{|(\sqrt{x(\xi)})'|}{\sqrt{x(\xi)}} + \sup\limits_{\xi \in [t_1;t]}\frac{|(\sqrt{x(\xi)})'|}{\sqrt{x(\xi)}}\le
$$
$$
\le \frac{c[1 + \sqrt{x(t_0)}(t_1 - t_0)]}{2[1 + \sqrt{x(t_0)}(t - t_0)]^{1 +2\alpha}} + \frac{c}{2[1 + \sqrt{x(t_0)}(t_1 - t_0)]^\alpha}, \phantom{aaa}t > \overline{t}.   \eqno (3.7)
$$
From definition of  $t_1(t)$ it follows, that  $t - t_1 < \frac{1}{2}(t - t_0)$  for  $t > \overline{t}$.
Then  $t_1 - t_0 > \frac{1}{2}(t - t_0)$    for  $t > \overline{t}$, and therefore from (3.7) we obtain:
$$
R_x(t_1(t);t)\le \frac{c}{2}\Biggl\{\frac{1}{[1 + \sqrt{x(t_0)}(t - t_0)]^\alpha} + \frac{1}{\bigl[1 +\frac{\sqrt{x(t_0)}}{2}(t - t_0)\bigr]^\alpha}\Biggr\}\le \frac{2^{\alpha - 1} c}{[1 + \sqrt{x(t_0)}(t - t_0)]^\alpha},
$$
 $t > \overline{t}$.  From here we immediately get:

{\bf Lemma 3.2.} {\it Let  $x(t)$ satisfies the conditions  (3.6). Then
$$
\rho_x(t) \le \frac{2^{\alpha - 1} c}{[1 + \sqrt{x(t_0)}(t - t_0)]^\alpha},\phantom{aaa} t > \overline{t}.\phantom{aaa}
\Box$$
}

Consider the sets
$$
A_t = A_t(x)\equiv \{s\in [t_0;t]: y'_x(s) \ge 0\},\phantom{aaa} B_t=B_t(x)\equiv \{s\in [t_0;t]: y'_x(s) < 0\}.
$$
It is evident, that  $A_t$  and  $B_t$  are measurable  and
$$
A_t \cup B_t = [t_0; t], \phantom{aaa} A_t \cap B_t =  \emptyset.  \eqno (3.8)
$$
Suppose $s\in A_t$. Then  $y'_x(s) \ge 0,\phantom{a} y_x(s) \le \sqrt{x(s)}$, and therefore
$$
\int\limits_{A_t}\frac{y'_x(s)}{2\sqrt{x(s)}}ds \le \int\limits_{A_t}\frac{y'_x(s)}{y_x(s) + \sqrt{x(s)}}ds\le \int\limits_{A_t}\frac{y'_x(s)}{2y_x(s)}ds. \eqno (3.9)
$$
For  $s\in B_t$ we have: $y'_x(s) < 0,\phantom{a} y_x(s) > \sqrt{x(s)}$. Then
$$
\int\limits_{B_t}\frac{y'_x(s)}{2\sqrt{x(s)}}ds \le \int\limits_{B_t}\frac{y'_x(s)}{y_x(s) + \sqrt{x(s)}}ds\le \int\limits_{B_t}\frac{y'_x(s)}{2y_x(s)}ds.
$$
Summarizing each  part of these inequalities with the corresponding parts of (3.9) and taking into account (3.8) we get:
$$
\int\limits_{t_0}^t\frac{y'_x(s)}{2\sqrt{x(s)}}ds \le \int\limits_{t_0}^t\frac{y'_x(s)}{y_x(s) + \sqrt{x(s)}}ds\le \int\limits_{t_0}^t\frac{y'_x(s)}{2y_x(s)}ds, \phantom{aaa}t\ge t_0.
$$
Due to the equality $y'_x(s) = (\sqrt{x(s)} - y_x(s))(\sqrt{x(s)} - y_x(s)),\phantom{a} s\ge  t_0$, from here we obtain the inequality
$$
-\frac{1}{2}\ln\frac{y_x(t)}{y_x(t_0)} \le \int\limits_{t_0}^t[y_x(s) - \sqrt{x(s)}] d s \le -  \int\limits_{t_0}^t\frac{y'_x(s)}{2\sqrt{x(s)}} d s, \phantom{aaa}t\ge t_0.
$$
Therefore,
$$
\frac{1}{4}\ln\left[\frac{x(t)}{y_x^2(t)}\right] \le \int\limits_{t_0}^t[y_x(s) - \sqrt{x(s)}]d s + \frac{1}{4}\ln \left[\frac{x(t)}{x(t_0)}\right] \le \int\limits_{t_0}^t \frac{(\sqrt{x(s)} - y_x(s))'}{2\sqrt{x(s)}} d s, \phantom{aaa}t\ge t_0.
$$
Then integrating the last integral by parts we obtain:
$$
\frac{1}{4}\ln \left[\frac{x(t)}{y^2_x(t)}\right] \le \int\limits_{t_0}^t[y_x(s) - \sqrt{x(s)}]d s + \frac{1}{4}\ln \left[\frac{x(t)}{x(t_0)}\right] \le \frac{1}{2} -
$$
$$
- \frac{y_x(t)}{2\sqrt{x(t)}} + \int\limits_{t_0}^t\frac{[\sqrt{x(s)} - y_x(s)]x'(s)}{4 x(s)^{3\slash 2}} d s.\phantom{aaa} t\ge t_0, \eqno (3.10)
$$
or
$$
\frac{1}{4}\ln \left[\frac{x(t)}{y^2_x(t)}\right] \le \int\limits_{t_0}^t[y_x(s) - \sqrt{x(s)}]d s + \frac{1}{4}\ln \left[\frac{x(t)}{x(t_0)}\right] \le \frac{1}{2} -
$$
$$
- \frac{y_x(t)}{2\sqrt{x(t)}} + \int\limits_{t_0}^t\frac{y'_x(s)x'(s)}{4[\sqrt{x(s)} + y_x(s)] x(s)^{3\slash 2}} d s. \phantom{aaa}t\ge t_0, \eqno (3.11)
$$
Consider the function
$$
Q_x(t) \equiv \int\limits_{t_0}^t[y_x(s) - \sqrt{x(s)}]d s + \frac{1}{4}\ln x(t),\phantom{aaa} t\ge t_0.
$$

 {\bf Lemma 3.3}. {\it Let  $x(t)$   be a monotone nondecreasing function, and let for some $\varepsilon > 0$    the function  $\frac{x'(t)}{x(t)^{3\slash 2 - \varepsilon}}$   be bounded. Then  $Q_x(t)$ is bounded.}

Proof.  Since  $x(t)$  is a monotone nondecreasing function, then (see  [13]) $y_x(t) \le \sqrt{x(t)}$. from here and from the first inequality of (3.11) it follows, that  $Q_x(t) \ge \frac{1}{4} \ln x(t_0) >\linebreak >  - \infty,\phantom{a} t\ge t_0$.  Therefore,  $Q_x(t)$   is bounded from below.
Suppose $\frac{|x'(t)|}{x(t)}^{3\slash 2 - \varepsilon} \le c,\phantom{a} t\ge t_0$,
for some  $\varepsilon > 0, \phantom{a} c > 0$.  Then taking into account the inequality  $y_x(t) \le \sqrt{x(t)},\phantom{a} t\ge t_0$,  we will have:
$$
\int\limits_{t_0}^t\frac{y'_x(s)x'(s)}{4[\sqrt{x(s)} + y_x(s)]x(s)^{3\slash 2}}d s \le
 \int\limits_{t_0}^t\frac{y'_x(s)x'(s)}{8y_x(s)^{1 + 2\varepsilon} x(s)^{3\slash 2 - \varepsilon}} d s \le \frac{c}{8} \int\limits_{t_0}^{+\infty}\frac{d(y_x(s))}{y_x(s)}^{1 + 2\varepsilon} \stackrel{def}= d_0 < +\infty,
$$
$ t\ge t_0.$
From here and from the second inequality of (3.11) it follows, that $Q_x(t) \le d_0 + \frac{1}{2} + \frac{1}{4} \ln x(t_0) < +\infty,\phantom{a} t\ge t_0$. Therefore  $Q_x(t)$ is bounded  above. The lemma is proved.

{\bf Lemma 3.4.} {\it Let  $x(t) \ge \varepsilon > 0,\phantom{a} t\ge t_0,$ $\frac{x'(t)}{x(t)}$  be bounded, and let  $\int\limits_{t_0}^{+\infty}\rho_x(s)\frac{|x'(s)|}{x(s)^{3\slash 2}}d s <\\ < +\infty$. Then $Q_x(t)$   is bounded.}

Proof. By virtue of mean value theorem
$$
\ln\left[\frac{x(t)}{y_x(t)}\right] = 2[\ln \sqrt{x(t)} - \ln y_x(t)] = 2 \frac{\sqrt{x(t)} - y_x(t)}{\xi(t)}, \eqno (3.12)
$$
where  $\xi(t) \in [\min\{\sqrt{x(t)}, y_x(t)\};\max\{\sqrt{x(t)}, y_x(t)\}], \phantom{a}t\ge t_0$.
Since $x(t) \ge \varepsilon ,\phantom{a} t\ge t_0$, we have   $\xi(t) \ge \min\{\sqrt{x(t)}, \phantom{a} y_x(t)\}\ge \sqrt{\varepsilon},\phantom{a} t\ge t_0$. From here, from the boundedness of $\frac{x'(t)}{x(t)}$ and from (3.4),
(3.5), (3.12) it follows:
$$
\left|\ln\left[\frac{x(t)}{y_x^2(t)}\right]\right| \le \frac{2}{\sqrt{\varepsilon}}\rho_x(t) \le \frac{c}{\sqrt{\varepsilon}} < + \infty,  \phantom{a} \mbox{where}\phantom{a} \frac{|x'(t)|}{x(t)} \le 2 c, \phantom{aaa} t\ge t_0.
$$
From here and from the first inequality of (3.10) it follows, that  $Q_x(t) \ge \frac{1}{4}\ln x(t_0) - \frac{c}{\varepsilon} > \linebreak > - \infty,\phantom{a} t\ge t_0$. Therefore,  $Q_x(t)$
is bounded  below. From (3.4) it follows:
$$
\int\limits_{t_0}^t\frac{[\sqrt{x(s)} - y_x(s)]x'(s)}{4x(s)^{3\slash 2}} d s \le \int\limits_{t_0}^{+\infty}\frac{\rho_x(s) |x'(s)|}{4 x(s)^{3\slash 2}} d s \stackrel{def}= d_1 < +\infty.
$$
Then taking into account the second inequality of (3.11) we will have: $Q_x(t) \le \frac{1}{2} + d_1 + \\ + \frac{1}{4}\ln x(t_0) < + \infty$. Therefore, $Q_x(t)$  is bounded  above. The lemma is proved.

Consider the Riccati equation
$$
y'(t) + y^2(t) = \frac{D(t)}{4},\phantom{aaa} t\ge t_0.   \eqno (3.13)
$$
In the sequel we shall assume, that the conditions A) are satisfied.  Each solution of Eq. (2.14), existing on $[t_0;+\infty)$,   is connected with some solution $\psi(t)$ of Eq.
 (1.3) by the equality (see [3], pp. 391, 392).
$$
\psi(t) = \psi(t_0)\exp\biggl\{\int\limits_{t_0}^ty(\tau)d\tau\biggr\},\phantom{aaa} t\ge t_0, \phantom{a} \psi(t_0) \ne 0.
$$
By  (1.2) from here it follows, that
$$
\phi_0(t) \equiv \exp\biggl\{\int\limits_{t_0}^t[y_{_{D\slash 4}}(\tau) - \frac{1}{2} p(\tau)] d\tau\biggr\}, \phantom{aaa}  t\ge t_0,
$$
is a solution of Eq. (1.1). Since  $\phi'_0(t) = [y_{_{D\slash 4}}(t) - \frac{1}{2} p(t)]\phi_0(t),\phantom{a} t\ge t_0$, we have
$$
|\phi_0'(t)| \le |y_{_{D\slash 4}}(t) - \frac{1}{2}\sqrt{D(t)}| |\phi_0(t)| + \frac{1}{2}|\sqrt{D(t)} - p(t)| |\phi_0(t)|. \phantom{aaa} t\ge t_0.
$$
By  (3.4) from here it follows:
$$
|\phi_0'(t)| \le \rho_{_{D\slash 4}}(t) |\phi_0(t)| + \frac{1}{2}\biggl[1 + |p(t) - \sqrt{D(t)}|\biggr] |\phi_0(t)|, \phantom{aaa} t\ge t_0. \eqno (3.14)
$$
Since
$$
\frac{1}{2}[\sqrt{D(t)} - p(t)] \phi_0(t) = \phi_0'(t) + \biggl[\frac{1}{2}\sqrt{D(t)} - y_{_{D\slash 4}}(t)\biggr] \phi_0(t),\phantom{aaa} t\ge t_0,
$$
we have
$$
\frac{1}{2}[1 + |p(t) - \sqrt{D(t)}|] |\phi_0(t)| \le |\phi_0'(t)| + [1 + |\sqrt{D(t)} - y_{_{D\slash 4}}(t)|] |\phi_0(t)|,\phantom{aaa} t\ge t_0.
$$
By virtue of (3.4) it follows from here, that
$$
[1 + |p(t) - \sqrt{D(t)}|] |\phi_0(t)| \le 2 |\phi_0'(t)| + [1 + \rho_{_{D\slash 4}}(t)] |\phi_0(t)|, \phantom{aaa} t\ge t_0.  \eqno (3.15)
$$
It is  not difficult to see, that
$$
|\phi_0(t)| = 2 \exp\{Q_{_{D\slash4}}(t) + \frac{1}{2}r_1(t)\},\phantom{aaa} t\ge t_0. \eqno (3.16)
$$
From here and from (3.14) it follows:
$$
|\phi_0'(t)| \le \rho_{_{D\slash4}}(t)|\phi_0(t)| + \exp\{Q_{_{D\slash4}}(t) + \frac{1}{2} r_2(t)\}, \phantom{aaa} t\ge t_0. \eqno (3.17)
$$
It follows from (3.15), that $\exp\{Q_{_{D\slash4}}(t) + \frac{1}{2}r_2(t)\} = \exp\biggl\{\int\limits_{t_0}^t[y_{_{D\slash 4}}(\tau) - \frac{1}{2} Re p(\tau)]d\tau +\\ + \ln [1 + |p(t) - \sqrt{D(t)}|] - \frac{1}{4}\ln 4\biggr\}\le [1 + |p(t) - \sqrt{D(t)}|] |\phi_0(t)| \le 2|\phi_0'(t)| + [1 + \rho_{_{D\slash4}}(t)] |\phi_0(t)|$,  $ t\ge t_0$.
Therefore,
$$
\exp\{r_2(t)\} \le \exp\{ -2Q_{_{D\slash4}}(t)\}\bigl[2 |\phi'_0(t)| + [1 +\rho_{_{D\slash4}}(t)] |\phi_0(t)|\bigr],\phantom{aaa} t\ge t_0. \eqno (3.18)
$$

 {\bf Lemma  3.5}. {\it  All solutions of Eq. (1.1) are bounded (vanish on  $+\infty$)  if and only if
the function  $\phi_0(t)$ is bounded (vanishes on  $+ \infty$).}

See the proof in [13].

 {\bf Lemma 3.6}. {\it  Eq. (1.1) is  Liapunov stable (asymptotically) if and only if  $\phi_0(t)$  and   $\phi'_0(t)$ are bounded (vanish on   $+\infty$).}

See the proof in [13].

{\bf Proof of Theorem 2.1}. Since the conditions A) hold by virtue of  Lemma 3.3 if the conditions  B) are satisfied, then  the function  $Q_{_{D\slash4}}(t)$  is bounded. If the conditions C) are satisfied, then the boundedness of  $Q_{_{D\slash4}}(t)$ follows from  Lemma 3.4. Thus the satisfiability of either  B) or  C) ensures the boundedness  of $Q_{_{D\slash4}}(t)$. Then from (3.15) it follows, that the function  $\phi_0(t)$ is bounded
($\lim\limits_{t\to +\infty}\phi_0(t) = 0$), if and only if the function $r_1(t)$  is bounded  above ($\lim\limits_{t\to+\infty} r_1(t)=-\infty$). By virtue of Lemma 3.5 from here it follows, that all solutions of Eq. (1.1) are bounded (vanish on  $+\infty$)
if and only if the function $r_1(t)$  is bounded  above  ($\lim\limits_{t\to+\infty} r_1(t)=-\infty$).   The theorem is proved.

{\bf Proof of Theorem 2.2}. From  D) follows  B). Therefore, by already proven the conditions A),  C) and  D) provide the boundedness of the function $Q_{_{D\slash4}}(t)$.  Then from (3.16) - (3.18) it follows, that  $\phi_0(t)$  and  $\phi'_0(t)$   are bounded (vanish on $+\infty$) if and only if $r_1(t)$  and  $r_2(t)$ are bounded  above ($\lim\limits_{t\to+\infty} r_j(t)= - \infty,\phantom{a} j=1,2$).
Since\linebreak  $r_1(t) \le r_2(t), t\ge t_0$,   from the  boundedness  above of $r_2(t)$
(from the equality $\lim\limits_{t\to+\infty} r_2(t)= - \infty$)
it follows the boundedness  above of  $r_1(t)$   (the equality\linebreak   $\lim\limits_{t\to+\infty} r_1(t)= - \infty$).  By virtue of  Lemma 3.6 from here it follows that Eq. (1.1) is  Liapunov stable (asymptotically) if and only if the function   $r_2(t)$ is bounded  above  ($\lim\limits_{t\to+\infty} r_2(t)=~- \infty$).
The theorem is proved.

{\bf Proof of Corollary 2.1}.  By virtue of  Lemma 3.2 from the first two conditions of corollary it follows
$$
\rho_{_{D\slash 4}}(t) \le \frac{c_1}{(1 + t - t_0)^\alpha},\phantom{a} t\ge t_0, \phantom{a}c_1 = const.
$$
Then
$$
\int\limits_{t_0}^{+\infty}\rho_{_{D\slash 4}}(\tau)\frac{|D'(\tau)|}{D(\tau)^{3\slash 2}} d \tau \le c_1\int\limits_{t_0}^{+\infty}\frac{d\tau}{\sqrt{D(s)}(1 + \tau - t_0)^{2\alpha}} < +\infty.
$$
Thus the group of conditions  C) of Theorem 2.1 is satisfied. Then А$_1$)
follows from  Theorem 3.1,  and   B$_1$) follows from  Theorem 2.2. The corollary is proved.

\pagebreak




\vskip 12pt
{\bf \centerline {References}}

\vskip 20pt

1. L. Cesary.  Asymptotic  behavior  and stability problems in  ordinary differential\\ \phantom{aaaaa} equations. Moskow, ''Mir'', 1964.

2. R. Bellman. Stability theory of differential equations, Moscow, Izdatelstvo inostrannoj \phantom{aaaaa} literatury, 1954.

3. Ph. Hartman.  Ordinary differential equations. Moscow, ''Mir'', 1970.

4. L. Ya. Adrianova. Introduction in the theory of linear systems of differential equations. \\ \phantom{aaaaa} St. Peterburg, Izdatelstvo St. Peterburgskogo universiteta, 1992.

5. M. V. Fedoriuk. Asymptotic methods for linear ordinary differential equations. \\ \phantom{aaaaa} Moskow, ''Nauka'', 1983.

6. N V. McLachlan. Theory and application of Mathieu functions. Moskow, ''Mir'',\\ \phantom{aaaaa} 1953.

7. I. M. Sobol. Study of the asymptotic behaviour of the solutions of the linear\\ \phantom{aaaaa} second order differential equations wit the aid of polar coordinates. "Matematicheskij \\ \phantom{aaaaa} sbornik", vol. 28 (70), N$^\circ$ 3, 1951, pp. 707 - 714.

8. L. A. Gusarov. Convergence to zero of solutions of linear second order differential \\ \phantom{aaaaa} equations. DAN SSSR, vol. LXXI,  N$^\circ$ 1, 1950, pp. 9 - 12.

9. G. A. Grigorian. Some properties of the solutions of linear second order ordinary\\ \phantom{aaaaa}  differential equations. Trudy UrO RAN, vol. 19, $N^\circ$ 19, 2013, pp. 69 - 90.

10.  G. A. Grigorian. Boundedness and stability criteria for linear ordinary differential\\ \phantom{aaaaa} equations of the second order. "Izvestia vuzov, Matematika, $N^\circ$ 12, 2013,\\ \phantom{aaaaa} pp. 11 - 18.

11. I. Knovles. On stability Conditions for Second Order Linear Differential Equations,\linebreak \phantom{aaaaa} Journal od Differential Equations 34, 179 - 203 (1979).

12. L. H. Erbe. Stability results for Periodic second Order Linear Differential Equations.\linebreak \phantom{aaaaa} Proc. Amer. Math. Soc., Vol. 93, Num. 2, 1985. pp. 272 - 276.

13. G. A. Grigorian. Some properties of differential root and theirs applications. Acta \linebreak \phantom{aaaaa} Math. Univ. Comenianae, Vol. LXXXV, 2 (2016), pp. 205 - 212.

14. G. A. Grigorian.  On two comparison tests for second-order linear  ordinary differential\linebreak \phantom{aaaaa} \hskip 2pt equations (Russian) Differ. Uravn. 47 (2011), no. 9, 1225 - 1240; translation in \linebreak \phantom{aaaaa} Differ. Equ. 47 (2011), no. 9 1237 - 1252, 34C10.

\end{document}